\documentclass[11pt]{amsart}
\usepackage{amsfonts,latexsym,amsthm,amssymb,graphicx}
\usepackage[all]{xy}
\usepackage{hyperref}
\usepackage[usenames]{color} 

\DeclareFontFamily{OT1}{rsfs}{}
\DeclareFontShape{OT1}{rsfs}{n}{it}{<-> rsfs10}{}
\DeclareMathAlphabet{\mathscr}{OT1}{rsfs}{n}{it}

\CompileMatrices

\newtheorem{theorem}{Theorem}[section]
\newtheorem{lemma}[theorem]{Lemma}

\newtheorem{prop}[theorem]{Proposition}
\newtheorem{corollary}[theorem]{Corollary}

{\theoremstyle{defin} \newtheorem{definition}[theorem]{Definition}}
{\theoremstyle{remark} \newtheorem{remark}[theorem]{Remark}
}

\numberwithin{equation}{section}

\title{Galois Extensions of Local Fields with One Wild Ramification Jump}
\author{Samuel Goodman}
\date{November 2025}
\begin{document}
\maketitle \vspace{-2.5 em} \begin{abstract}
For a given positive integer $n$ and $K/\mathbb{Q}_p$ a finite extension of ramification degree $e$, we determine the number of finite Galois extensions $L/K$ with inertia degree $f$ and a single nonnegative ramification jump at $n$ as long as $(p,e)$ is outside of a finite set. This builds upon the tamely ramified case, which is a classical consequence of Serre's Mass Formula, exhibiting a more restrictive behavior than in the tamely ramified case because the degrees of such extensions are bounded. We do this by working in a fixed Lubin-Tate extension and exploiting the surjectivity of a map corresponding to the ramification jump to reconstruct the $U^1$ part of the norm subgroup (coming from local class field theory) from its fibers and then by understanding how the fibers interact by studying them in terms of properties of the formal logarithm and partitions. 
\end{abstract}

\section{Introduction}
Fix a local field $K$, for example a finite extension of $\mathbb{Q}_p$. A classical question is counting the number of tamely ramified extensions of $K$ of a given degree, which in turn is done by Serre's mass formula (see [5]), giving precisely $n$ degree $n$ totally tamely ramified extensions (where $p \nmid n$), and then since a tamely ramified extension is uniquely a totally tamely ramified extension of its maximal unramified one and there is a unique unramified degree $k$ extension of $K$ for each $k$, we find there are $\sigma_0(n)$ tamely ramified degree $n$ extensions, where $\sigma_0(n)$ is the sum of the divisors of $n$ coprime to $p$. \newline \newline Using Kummer theory and the explicit description of tamely ramified Galois extensions or more elementary techniques (see [2]), one can also obtain the following result: 
\newline \newline Let $K/\mathbb{Q}_p$ be a finite extension and let $L$ be a finite unramified extension of $K$. Set $q=|k|$. Then there are $\gcd(n,q-1)$ tamely ramified Galois extensions $M/K$ such that $L \subset M$ and $M/L$ is totally ramified of degree $n$. \newline \newline This Kummer Theory approach quickly devolves into a certain counting of fixed points under the Galois group of the unramified part, something which we greatly build upon using class field theory in a broader context. \newline \newline Recall the definition of ramification groups: given a Galois extension $L/K$ of local fields with Galois group $G$ such that $L$ has ring of integers $\mathcal{O}_L$ and maximal ideal $\mathfrak{m}$, we define $G_i=\{\sigma \in G, \forall \alpha \in \mathcal{O}_L, \sigma(\alpha) \equiv \alpha \mod \mathfrak{m}^{i+1}\}$. We say that $L/K$ has a ramification jump at $i$ if $G_i \neq G_{i+1}$. In terms of ramification jumps, tamely ramified Galois extensions can be thought of as $(-1,0)$ Galois extensions since they have ramification jumps only at $-1$ and $0$. The next natural question to ask is how to count $(-1,n)$ extensions for a given positive integer $n$; these are the extensions with only a single wild ramification jump. \newline \newline In this paper, we completely resolve the case of one wild jump extensions for $n=2$ and resolve up to finitely many $(p,e)$ pairs, namely those with $p-1 \leq e<n$, for general $n$. In fact, we determine explicit formulas in the cases in all such cases. The general result is (see Theorem $7.1$):
\begin{theorem}
Suppose $K/\mathbb{Q}_p$ is finite with ramification index $e$, inertia degree $f'$, and residue field $k$. Let $f$ be a positive integer and let $q=p^{f'}, q'=q^f$, $k_0=\mathbb{F}_{q'}$, and $G=\text{Gal}(k_0/k)$. Then let $S$ be the set of $G$-invariant subspaces of $k_0$ and $t \in S$ to be the trace $0$ subspace under $Tr_{k_0/k}$. Let $b=n-1-\lfloor \frac{n-1}{p} \rfloor$. Then the number of $(-1,n)$ Galois extensions $L/K$ of inertia degree $f$ is $$\begin{cases} 0 & p|n, e \geq n \\ q^{bf'f+1}\sum_{h \in S} \frac{|h \cap t|}{|h|^{bf'+1}} & e \geq n, p \nmid n \\ 0 & e<n, p>e+1 \end{cases}$$ \end{theorem}
For each $n$, this theorem completely settles the question for all but finitely many $(p,e)$ pairs. In fact, for $n=2$, only one pair is left open, namely $e=1,p=2$, which is the unramified case at $p=2$. The result in this case is (see Theorem $7.3$): 
\begin{theorem}
Suppose $K/\mathbb{Q}_2$ is unramified of degree $n$. Then the number of $(-1,2)$ Galois extensions $L/K$ with inertia degree $r$ is $2^{n+1}$ (independent of $r$). \end{theorem} We notice an interesting dichotomy when $n=2$, which is that there is an obstruction to $2$-adic extensions as soon as the base extension $K/\mathbb{Q}_p$ becomes ramified while there is an obstruction everywhere else when $K/\mathbb{Q}_p$ is ramified (and no longer any obstruction at $p=2$). We also see that number of extensions grows rapidly based on the degree of the unramified part in the tower and that in general, the count, when nonzero, depends on the full data of the residue field extension in Theorem $1.1$ but suddenly only depends on the base residue field in Theorem $1.2$. \newline \newline Theorem $1.1$ also specializes to a particularly clean closed form when $K/\mathbb{Q}_p$ is totally ramified, giving the following result (see Corollary $7.2$): \begin{corollary} Suppose that $K/\mathbb{Q}_p$ is totally ramified. For a positive integer $f$, set $m(x)=\frac{x^f-1}{(x-1)^{p^{v_p(f)}}} \in \mathbb{F}_p[x]$, $\text{deg}(m)=f-p^{v_p(f)}=d$, and $\zeta_m(s)=\sum_{i=0}^{d} \frac{a_i}{p^{is}}$ be the $\zeta$-function for the ring $S=\frac{\mathbb{F}_p[x]}{(m)}$ (equivalently, $a_i$ is the number of degree $i$ monic factors of $m$). Then the number of $(-1,n)$ Galois extensions $L/K$ with inertia degree $f$ is $$\begin{cases} 0 & p|n, e \geq n \\ \frac{(p^{(f+1)b}-p^{fb}+p^{(f-1)b+1}-p^{db+1})}{p^b-1}\zeta_m(b)  & e \geq n, p \nmid n \\ 0 & e<n, p>e+1 \end{cases}$$ \end{corollary}
\noindent The intuition for the proof is as follows: Local Class Field Theory turns the problem into a counting problems on norm groups, and to track Galoisness, we just need to track the Galois action of the unramified part on the corresponding norm subgroup (see Proposition $2.2$). Then working within a fixed Lubin-Tate extension, the ramification jump conditions tell us that a projection map on the corresponding norm group is surjective in Proposition $3.1$. The crucial insight is to leverage this surjectivity by working with the fibers of the map itself and showing that the subgroup can be reconstructed from the fibers (as long as they satisfy certain compatibility conditions) and vice versa, as completed in Proposition $6.4$. These fibers can be expressed in a way that make them compatible with each other in the $e \geq n$ setup due to the isomorphism of Lemma $6.3$. The resulting compatibility relations are more combinatorial, allowing one to use formal power series approach combined with some combinatorics involving partitions (as in Section $4$). Combined with some structure theory of characteristic $p$ local fields, one can understand the behavior of such compatibility relations entirely, as in Theorem $5.1$. This construction relies crucially on the condition $e \geq n$ (and in fact, comparing to Theorem $7.2$, cannot possibly work in many cases when $e<n$). Finally, we patch together the information coming from within each individual Lubin-Tate extension by looking at a choice of uniformizer up to units and counting fixed points, which is ring theory, leading to the proof of Theorem $7.1$. Then most of the remaining cases are shown to be infeasible because the norm subgroups are forced to be too big. We then modify the construction in the unramified case at $p=2$, the only exceptional case when $n=2$, exploiting the fact that there is only one dimension of flexibility.
\section{Preliminaries}
We start by noting an immediate consequence of [4, IV, Prop. 7]:
\begin{prop} 
Let $L/K$ be a Galois extension with Galois group $G$. Set $G_0=I$ and let $G_i$ be the $i$th higher ramification group for $i \geq 1$. Then we have that: \newline \newline 
1) $G_0/G_1$ is isomorphic to a subgroup of $l^{\times}$, where $l$ is the residue field of $L$, and thus is cyclic. \newline \newline
2) For $i \geq 1$, $G_i/G_{i+1} \cong (\mathbb{Z}/p\mathbb{Z})^k$ for some $k$, and so in particular, $G_1$ is a $p$-group. 
\end{prop}
Letting $\pi$ be a uniformizer for $L/K$ and assuming the extension $L/K$ is totally ramified, we have that $\mathcal{O}_L=\mathcal{O}_K[\pi]$ and so the ramification groups are determined by $\pi$. More precisely, letting $\text{Gal}(L/K)=G$, we have that $G_n=\{\sigma \in G, \sigma(\pi) \equiv \pi \mod \mathfrak{m}^{n+1}\}$, where $\mathfrak{m}=(\pi)$ is the maximal ideal of $\mathcal{O}_L$. In general, we can also define $G_t=G_{\lceil t \rceil}$ so that $G_t$ is defined for all real $t \geq -1$. We will later seek to better understand these ramification groups. \begin{prop}
Let $M/L/K$ be a series of finite extensions of local fields such that $M/L$ is abelian and $L/K$ is Galois. Let $H$ be the norm subgroup of $L^\times$ associated to $M$ under Local Class Field Theory. Then $M/K$ is Galois iff $\sigma(H)=H$ for all $\sigma \in \text{Gal}(L/K)$.
\begin{proof} $M/K$ being Galois is equivalent to $\sigma(M)=M$ for all $\sigma: M \rightarrow \bar{K}$ an embedding fixing $K$. Since $L/K$ is Galois, $\sigma(L)=L$ for all such embeddings and thus $\sigma(M)/L$ is an abelian extension. Its norm group is clearly $\sigma(H)$. However, Local Class Field Theory gives an order-reversing bijection between norm groups and finite abelian extensions, and so we have that $\sigma(M)=M$ for all $\sigma$ iff $\sigma(H)=H$ for all $\sigma$. But $\sigma|_L$ precisely attains the elements of $\text{Gal}(L/K)$, and so we conclude. \end{proof} \end{prop} We will also need some preliminaries on the structure of multiplicative groups of local fields. Firstly, to understand how $p$th powers behave in local fields, we have (from [1, I, Prop. 6.2]:
\begin{prop} Let $K$ be a nonarchimedean local field of characteristic $p$ with ring of integers $\mathcal{O}_K$ and maximal ideal $\mathfrak{m}$, and fix a uniformizer $\pi$. Let $e_1,\cdots,e_f$ be an $\mathbb{F}_p$-basis for its residue field. Then the group $1+\mathfrak{m}$ has a topological basis (as a $\mathbb{Z}_p$-module) of the form $1+e_i \pi^j$ where $1 \leq i \leq f$ and $p \nmid j \in \mathbb{N}$. \end{prop} In one of our ranges, where $e \geq n$, a connection to the characteristic $p$ setting, where there is the above particularly clean topological basis, will give us quick control over certain Galois-equivariant homomorphisms later on. \newline \newline In our other range, namely when $p>e+1, e<n$, it will be useful to understand how $p$th powers behave. From [3, II, Proposition 5.5], we have that: \begin{prop} Let $K/\mathbb{Q}_p$ be finite with ring of integers $\mathcal{O}_K$, maximal ideal $\mathfrak{m}$ and ramification index $e$ and suppose that $p>e+1$. Then the logarithm map $1+\mathfrak{m} \to \mathfrak{m}$ is an isomorphism. \end{prop}
\begin{remark} According to this isomorphism, the $p$th powers $(1+\mathfrak{m})^p$ will correspond to $p\mathfrak{m}=\mathfrak{m}^{e+1}$, which will then correspond back to $1+\mathfrak{m}^{e+1}$, and so we have that $(1+\mathfrak{m})^p=1+\mathfrak{m}^{e+1}$ in this setting. \end{remark}
\section{Ramification Groups}
We now prove some general results about ramification groups. Let $\phi_{L/K}(t)=\int_{0}^{t} \frac{1}{[G_0:G_t]}dt$. This is known as the Hasse-Herbrand function and is an increasing piecewise linear $[-1,\infty) \to [-1, \infty)$ and is in particular bijective. We can then set $\psi_{L/K}$ to be the inverse function, which is also increasing and piecewise linear. Then we have Herbrand's Theorem, which states that $G_uH/H=(G/H)_v$, where $v=\psi_{L/K}(u)$. We are interested in the case of Galois extensions where there is a single jump in the wild ramification groups. We now characterize such extensions in the context of Lubin-Tate Theory. As before, we have that $G=\text{Gal}(K_{\pi,n}/K) \cong \mathcal{O}^{\times}_K/(1+\mathfrak{m}^n)$ and so view subgroups of $G$ in terms of subgroups of $\mathcal{O}^{\times}_K/(1+\mathfrak{m}^n)$. We see that an subextension $L/K$ with corresponding subgroup $H$ has a (lower) ramification jump at $u$ if $(G/H)_u \neq (G/H)_{u+1}$.
\begin{prop}
Suppose that $K_{\pi,n}/K$ is a Lubin-Tate extension with Galois group $G$. Assuming $n \geq k+1$, the subextensions $L/K$ with a single wild ramification jump at $k$ correspond to the proper subgroups $H$ of $G$ that both contain $(1+\mathfrak{m}^{k+1})/(1+\mathfrak{m}^n)$ and have the property that the canonical map $H \rightarrow \mathcal{O}_K^{\times}/(1+\mathfrak{m}^k)$ is surjective under the identification $G \cong \mathcal{O}_K^{\times}/(1+\mathfrak{m}^n)$. \end{prop}
\begin{proof} First recall from Lubin-Tate theory that if $m<q^n, q^r \leq m<q^{r+1}$, we have that $G_m=(1+\mathfrak{m}^r)/(1+\mathfrak{m}^n)$. It follows that $G^{r+1}=G_{q^r}$ for each $0 \leq r<n$. Next note that Galois extensions with a single positive ramification jump have $G_{\phi_{L/K}(v)}=G_v$ for all $v$. Indeed, letting the jump be at $k$, we have that for $v \leq k$, $G_{\phi_{L/K}(v)}=G_v$, where $\phi_{L/K}(v)=v$ follows since $[G_0:G_v]=1$ for all $0 \leq v<k+1$. On the other hand, for $v>k+1$, both are trivial since then $k<\phi_{L/K}(v)<v$. Thus the upper numbering and lower numbering groups coincide for such extensions.

\vspace{0.6\baselineskip} 
\noindent There being a unique jump at $k$ is equivalent to $(G/H)^{k+1}=1$ and $(G/H)^{k'}=(G/H)^0$ for $1 \leq k' \leq k$ by definition as these coincide with the lower ramification groups. By Herbrand's Theorem, this is equivalent to having $G^{k+1}=G_{q^k} \subset H$ and $G_{q^{k'}}H/H=(G/H)_0$ for $1 \leq k' \leq k$. The former condition is equivalent to containing $(1+\mathfrak{m}^{k+1})/(1+\mathfrak{m}^n)$. The latter just means that $G_{q^{k'}}H$ is constant for $1 \leq k' \leq k$, which is equivalent to having the condition that $H \rightarrow \mathcal{O}_K^{\times}/(1+\mathfrak{m}^{k'})$ is surjective for $1 \leq k' \leq k$. However, surjectivity at $k$ implies surjectivity elsewhere, implying the claim. \end{proof} 
\noindent We now use Lemma $3.1$ to say more about the $H$ such that there is a single jump at $1$. For $1$, the second condition is superfluous, and so it is enough to contain $(1+\mathfrak{m}^2)/(1+\mathfrak{m}^n)$. By Proposition $3.1$ and class field theory, this implies that $L \subset K_{\pi,2}$, and so these are precisely the working extensions. If we assume that there only a single jump at $1$, then we must also have that the map $H \rightarrow \mathcal{O}^\times_K/(1+\mathfrak{m})$ is surjective. However, we know that $H$ is isomorphic to a subgroup of $\mathcal{O}^\times_K/(1+\mathfrak{m}^2)$, which has $p$-Sylow group $(1+\mathfrak{m})/(1+\mathfrak{m}^2)$ and cyclic subgroup of order $q-1$ generated by the coset of $\mu_{q-1}$, giving a splitting $\mathcal{O}^\times_K/(1+\mathfrak{m}^2) \cong k^+ \times k^\times$, where $k$ is the residue field of $K$, where the isomorphism is induced upon fixing a uniformizer $\pi$. Any subgroup of this group will be isomorphic to the direct product of its Sylow subgroups, thus isomorphic to a product of subgroups $k_1,k_2$ of $k^+, k^\times$, respectively. \vspace{0.6\baselineskip} 

\noindent Surjectivity is then equivalent to $k_2$ being all of $k^\times$. Now let $a \in k$ and note that the coset $H_a$ of $H$ of elements congruent to $a \mod \mathfrak{m}$ is just $a+f(a)\pi \mod \pi^2$ for some unique coset $f(a)$ of $k^+/k_1$. We must then have that $H_aH_b=H_{ab}$, which implies that $af(b)+bf(a) \equiv f(ab) \mod k_1$. Now since there a unique subgroup $H$ with a given $k_1$ and every $f$ gives a different subgroup with that corresponding $k_1$, we conclude that for any subgroup $h$ of $k^+$, the unique ``differential" of the form $f: k \rightarrow k^+/k_1$ is the zero differential. Thus the $1$ case is encapsulating the differential information of the residue field. The case of $n$ can thus be seen as a more complicated type of differential. \newline \newline Now we investigate this case. First we need some lemmas on formal logarithms:
\section{Formal Logarithms}
\begin{lemma} Let $K$ be a characteristic $0$ field and consider the formal power series of $1+\sum_{i=1}^{\infty} a_iT^i \in K[T,a_1,\cdots,a_n, \cdots]$ in indeterminates $a_i$ and $T$. For a partition $p_j$ of $j$, let $\ell(p_j)$ be the length of the partition, let $p_{ji}$ be the number of occurrences of $i$ in the partition, and let $m(p_j)$ be the number of permutations of the partition. Let $P_j$ be the set of partitions of a positive integer $j$ ($P_{j,k}$ for those with length $k$) and for $p_j$ a partition of $j$. Then we have that $\log(1+\sum_{i=1}^{\infty} a_iT^i)=\sum_{n=1}^{\infty} (\sum_{p_n \in P_n} (-1)^{\ell(p_n)} \frac{m(p_n)}{\ell(p_n)}\prod_{i=1}^{n} a_i^{p_{ni}})T^n$. \end{lemma}
\begin{proof} To start, note that $$\log(1+\sum_{i=1}^{\infty} a_iT^i)=\sum_{m=1}^{\infty} \frac{(-1)^m}{m}(\sum_{e_1,\cdots,e_m} (\prod_{i=1}^{m} a_{e_i})T^{e_1+\cdots+e_m})$$ Fixing an exponent $e_1+\cdots+e_m=n$, the formal sum rearranges as $$\sum_{n=1}^{\infty} (\sum_{k=1}^{n} \frac{(-1)^k}{k} \sum_{e_1+\cdots+e_k=n} \prod_{m=1}^{k} a_{e_m})T^n$$ Now letting $f_i$ be the number of occurrences of $i$ among the $e_m$s, the product can be written $\prod_{m=1}^{j} a_i^{f_i}$ for each tuple $e_1,\cdots,e_k$. The exponents $f_i$ satisfy $\sum_{i=1}^{n} f_i=k$ and $\sum_{i=1}^{k} if_i=n$, giving a corresponding partition $p$ of $j$ with $f_i$ occurrences of $i$ so that $\ell(p)=k$ and $m(p)$ the number of ways in which the product can occur. Thus breaking down $\sum_{e_1+\cdots+e_k=n} \prod_{m=1}^{k} a_{e_m}$ by partitions gives $$\sum_{p_n \in P_{n,k}} m(p_n)a_i^{p_{ni}}$$ Finally, summing over $k$, corresponding to partition lengths, gives the desired coefficient for $T^j$. \end{proof} Note that the coefficient  $\sum_{p_n \in P_n} (-1)^{\ell(p_j)} \frac{m(p_n)}{\ell(p_n)}\prod_{i=1}^{j} a_i^{p_{ni}}$ is a multivariate polynomial in precisely $a_1,\cdots,a_n$. We wish to better understand these polynomials, specifically their integrality properties. Let $r_n(x_1,\cdots,x_n)=\sum_{p_n \in P_n} (-1)^{\ell(p_n)} \frac{m(p_n)}{\ell(p_n)}\prod_{i=1}^{n} x_i^{p_{ni}}$ and for all $n$ \begin{lemma} We have that $r_n(x_1,\cdots,x_n) \in \mathbb{Z}[\frac{1}{n}][a_1,\cdots,a_n]$ \end{lemma} \begin{proof} It suffices to show that for each partition $p_n \in P_n$, the rational number $\frac{m(p_n)}{\ell(p_n)}$ cannot have reduced denominator divisible by any $q \nmid n$. Choose some $q|\ell(p_n)$. We prove that $q^{v_q(\ell(p_n)}|m(p_n)$. Indeed, consider the set $T_{p_n}$ of all distinct permutations of the partition $p_n$ and consider the action on $T_{p_n}$ of the subgroup generated by $(1 \cdots \ell(p_n)) \in S_{\ell(p_n)}$, the symmetric group on $\ell(p_n)$ elements. For any fixed $t \in T_{p_n}$, the Orbit-Stabilizer Theorem implies that $|\text{Orb}(t)||\text{Stab}(t)|=\ell(p_n)$. Now suppose that $q||\text{Stab}(t)|$. Then there would be an element of order $q$ in $\text{Stab}(t)$, which would break the permutation $t$ into a block repeated $q$ times. But then the sum of the elements in $t$, which by definition is $j$ would be a multiple of $q$, contradicting $q \nmid j$. Thus $q$ cannot divide $|\text{Stab}(t)|$ and so $q^{v_q(\ell(p_n)}||\text{Orb}(t)|$. But the set $T_{p_j}$ is a disjoint union of its orbits under this action, and thus $q^{v_q(\ell(p_n)}||T_{p_n}|=m(p_n)$, as desired. \end{proof} \begin{remark} Lemma $4.2$ can also be seen in terms of big Witt vectors and ghost components: indeed, the formal series $1+\sum_{i=1}^{\infty} a_iT^i$ can be written uniquely as $\prod_{i=1}^{\infty} (1-x_iT^i)^{-1}$, where each $x_n$ is an integer polynomial in $a_1,\cdots,a_n$, and then taking logarithms gives the $T^n$ coefficient is $\frac{1}{n}\sum_{d|n} dx_d^{\frac{n}{d}}=\frac{W_n}{n}$ (in terms of ghost components), which implies that if $\log(1+\sum_{i=1}^{\infty} a_iT^i)=\sum_{i=1}^{\infty} c_iT^i$, then $nc_n \in \mathbb{Z}[a_1,\cdots,a_n]$ for each $n$. However, the partition identity is more convenient for explicitly tracking the coefficients, for example in the proof that $r_n(x_1,\cdots,x_n)-x_n$ satisfies the functional equation of Theorem $5.1$. \end{remark} We will now see that the polynomials $r_n(x_1,\cdots,x_n)=\sum_{p_n \in P_n} (-1)^{\ell(p_n)} \frac{m(p_n)}{\ell(p_n)}\prod_{i=1}^{n} a_i^{p_{ni}}$ satisfy a remarkable property. \begin{lemma}  Let $z_i=x_i+y_i+\sum_{j=1}^{i-1} x_jy_{i-j}$. Then in the polynomial ring $\mathbb{Q}[x_1,\cdots,x_n,y_1,\cdots,y_n]$, we have the identity $r_n(x_1,\cdots,x_n)+r_n(y_1,\cdots,y_n)=r_n(z_1,\cdots,z_n)$ \end{lemma}
\begin{proof} First note by additivity of the formal logarithm that $$\log(1+\sum_{i=1}^{\infty} x_iT^i)+\log(1+\sum_{i=1}^{\infty} y_iT^i)=\log(1+\sum_{i=1}^{\infty} (x_i+y_i+\sum_{j=1}^{i-1} x_jy_{i-j})T_i)=\log(1+\sum_{i=1}^{\infty} z_iT^i)$$ The claim then follows from Lemma $4.1$. \end{proof}
\section{Finite Fields}
In this section, we will first prove a general theorem about a functional equation over our finite fields and then prove some specialized lemmas in the $n=2$ case at the prime $2$ which will be essential in completing the only case when $n=2$ not covered by the general theorems. Following the convention of the previous section, let $c_i=a_i+b_i+\sum_{j=1}^{i-1} a_jb_{i-j}$ and $r_n(x_1,\cdots,x_n)=\sum_{p_n \in P_n} (-1)^{\ell(p_n)} \frac{m(p_n)}{\ell(p_n)}\prod_{i=1}^{n} x_i^{p_{ni}}$.
\begin{theorem} Suppose we have an extension $k_0/k$ of finite fields with Galois group $G$ and a Galois invariant subgroup $h \subset k_0^+$ and set $q=|k|$. Then the number Galois equivariant functions $f: k_0^{n-1} \to k_0^+/h$ such that $$f(a_1,\cdots,a_{n-1})+f(b_1,\cdots,b_{n-1})+\sum_{i=1}^{n-1} a_ib_{n-i} \equiv f(c_1,\cdots,c_{n-1}) \mod h$$ is $$\begin{cases} 0 & p|n \\ q^{(n-1-\lfloor \frac{n-1}{p} \rfloor)\text{codim}_{\mathbb{F}_p}(h)} & p \nmid n  \end{cases}$$ \end{theorem}
\begin{proof} First suppose that $p \nmid n$. By Lemma $4.2$, the polynomial $f(x_1,\cdots,x_{n}) \in \mathbb{Q}[x_1,\cdots,x_{n}]$ is actually an element of $\mathbb{Z}[\frac{1}{n}][x_1,\cdots,x_{n}]$, and so since $p \nmid n$, we can view it as an element in $k[x_1,\cdots,x_n]$. Then considering $r_n(x_1,\cdots,x_{n})-x_n$ and invoking Lemma $4.1$, removing the term corresponding to the partition $\{n\}$ of $n$, corresponding to the unique partition of $n$ containing $n$, the resulting polynomial is now a function of $x_1,\cdots,x_{n-1}$, and so set $g(x_1,\cdots,x_{n-1})=r_n(x_1,\cdots,x_n)-x_{n}$. By Lemma $4.4$, we have that $$r_n(x_1+\cdots,x_{n})+r_n(y_1+\cdots,y_{n})=r_n(z_1,\cdots,z_{n}) \in \mathbb{Q}[x]$$ It follows that $$g(x_1,\cdots,x_{n-1})+x_{n}+g(y_1,\cdots,y_{n-1})+y_{n}=g(z_1,\cdots,z_{n-1})+z_{n}$$ Since $z_{n}=x_{n}+y_{n}+\sum_{i=1}^{n-1} x_iy_{n-i}$, we see that $-g$ provides a solution. Furthermore, note that the reduction $-g$ lives in $\mathbb{F}_p[x_1,\cdots,x_{n-1}]$ and so it automatically gives a Galois-equivariant solution. Then any other solution $f_1$ can be written uniquely as $-g+f_2$, where $f_2$ is a solution to $$f(x_1,\cdots,x_{n-1})+f(y_1,\cdots,y_{n-1}) \equiv f(z_1,\cdots,z_{n-1}) \mod h$$ and so the number of solutions is just the number of possible equivariant $f_2$. But Galois-equivariant solutions to $$f(x_1,\cdots,x_{n-1})+f(y_1,\cdots,y_{n-1}) \equiv f(z_1,\cdots,z_{n-1}) \mod h$$ correspond to Galois-equivariant group homomorphisms $(1+Tk_0[T])/(1+T^nk_0[T]) \to k_0/h$ where the bijection comes from setting $f(x_1,\cdots,x_{n-1})=g(1+x_1T+\cdots+x_{n-1}T^{n-1})$, and so the number of solutions is $$|\text{Hom}_G((1+Tk_0[T])/(1+T^nk_0[T]), k/h)|$$ Then note that the group $k_0/h$ is $p$-torsion, so any such homomorphism will factor through the $p$-powers. However, choosing a Galois equivariant (using the Galois module structure of $k_0$) $\mathbb{F}_p$-basis for $k_0$, the topological generators of Proposition $2.3$ then give an $\mathbb{F}_p$-basis for $(1+Tk_0[T])/(1+Tk_0[T])^p$ of the form $$\{1+e_iT^j, 1 \leq i \leq f, p \nmid j \in \mathbb{N}\}$$ and thus quotienting by $1+T^nk_0[T]$ leaves precisely those with $j<n$. Now for fixed $j$ a Galois-equivariant homomorphism is determined on each Galois orbit on the elements $1+e_i T^j$ by a single element of the orbit, and there are $[k':\mathbb{F}_p]$ such orbits, for each of which there are $|k_0/h|$ choices for the image of that basis element. Finally, there are $n-1-\lfloor \frac{n-1}{p} \rfloor$ choices of $j$, giving $q^{(n-1-\lfloor \frac{n-1}{p} \rfloor)\text{codim}_{\mathbb{F}_p}(h)}$ total solutions, as desired. \newline \newline Now suppose that $p|n$. We will show that there exist no solutions. Suppose that $f$ is a solution. We will derive a contradiction by iterating a certain element with itself $p$ times. Indeed, start with some $(x_1,\cdots,x_{n-1})$ and set $x_{i,1}=x_i$. Then inductively set $x_{i,k}$ to be the $c_i$ determined by taking $a_i=x_{i,(k-1)}$ and $b_i=x_i$. Then we have that $$1+\sum_{i=1}^{n-1} x_{i,k}T^i=(1+\sum_{i=1}^{n-1} x_iT^i)^k \mod T^n$$ On the other hand, given the power series $1+\sum_{i=1}^{n-1} x_iT^i$, one can define new variables $x_{n,k}$ so that $x_{n,1}=0$ and in general $x_{n,k}$ is the $T^n$ coefficient of $(1+\sum_{i=1}^{n-1} x_iT^i)^k$. It then follows that the cross term in the functional equation for that level becomes $x_{n,k}-x_{n,1}-x_{n,(k-1)}$ because we are always doing the full multiplication but then removing the degree $n$ terms that cancel out, as in the proof of the other direction. Now choose $(a_1,\cdots,a_n)$ so that the first $\frac{n}{p}-1$ terms are $0$ and allow $a_{\frac{n}{p}}$ to vary over $k$. Iterating the functional equation $p$ times, the sum of all cross terms will be $$\sum_{i=2}^{p} x_{n,i}-x_{n,i-1}-x_{n,1}$$ which telescopes to $x_{n,p}-px_{n,1}=x_{n, p}$. We can compute $x_{n, p}$ as the $T^n$ coefficient in $$(1+a_{\frac{n}{p}}T^{\frac{n}{p}}+\sum_{i=\frac{n}{p}+1}^{n-1} a_iT^i)^p$$ By Freshman's Dream, we just get $1+a_{\frac{n}{p}}T^n+\sum_{i=\frac{n}{p}+1}^{n-1} a_i^pT^{ip}$ and so the only contribution to $T^n$ will be the term $a_{\frac{n}{p}}^{p}$, which can attain any value in $k$ as $a_{\frac{n}{p}}$ varies. On the other hand, writing the functional equation as $$\sum_{i=1}^{n-1} a_ib_{n-i} \equiv f(c_1,\cdots,c_{n-1})-f(a_1,\cdots,a_{n-1})-f(b_1,\cdots,b_{n-1}) \mod h$$ and summing the $f$ values also leads to telescoping, giving $$f(a_{1,p},\cdots,a_{n-1,p})-pf(a_1,\cdots,a_{n-1}) \equiv f(a_{1,p},\cdots,a_{n-1,p}) \equiv 0 \mod h$$ because there are no terms of degree below $n$. Thus whenever $h$ is a proper subgroup of $k_0^+$, we can choose $a_{\frac{n}{p}}$ to give a contradiction, completing the proof. \end{proof}
Now we prove a slightly refined version of Theorem $5.1$ but only when $n=2,p=2$.
\begin{lemma} Let $k_0$ be a finite field of characteristic $2$ and $h,h'$ subgroups of $k_0^+$. Let $r=\dim_{\mathbb{F}_p}(h')$. Say that $h \sim h'$ if $x^2 \in h$ for each $x \in h'$. Then the number of functions $f: h' \rightarrow k_0^+/h$ such that $f(a+b) \equiv f(a)+f(b)+ab \mod h$ is $$\begin{cases} 0 & h \not \sim h' \\ 2^{r(\text{codim}_{\mathbb{F}_p}(h))} & h \sim h'  \end{cases}$$ \end{lemma} \begin{proof} From this relation and an easy induction, we deduce $$f(\sum_{i=1}^{r} a_i) \equiv \sum_{i=1}^{r} f(a_i)+\sum_{1 \leq i<j \leq r} a_ia_j \mod h$$ which in particular implies that for any positive integer $r$, $f(ra) \equiv rf(a)+\binom{r}{2}a^2 \mod h$ upon setting all $a_i$s equal. If $h \sim h'$, we find that $f(2a) \equiv a^2 \equiv 0 \mod h$ and so $f$ is in fact well-defined. If $p=2$ and $h \not \sim h'$, then we get a contradiction since we would need $a^2 \equiv 0 \mod h$ for all $a \in h'$, meaning that no such functions can exist. \newline \newline Now let $e_1,\cdots, e_r$ be an $\mathbb{F}_2$-basis for $h'$. Upon choosing $f(e_i)$, the above relation gives $f(re_i) \equiv rf(e_i)+\binom{r}{2}e_i^2 \mod h$, and so $f(re_i)$ is determined by $f(e_i)$. Furthermore, we must have $$f(\sum_{i=1}^{n} a_ie_i) \equiv \sum_{i=1}^{n} a_if(e_i)+\binom{a_i}{2}e_i^2+\sum_{1 \leq i<j \leq n} a_ia_je_ie_j \mod h$$ and so $f$ is completely determined by $f$ on the basis. 
\newline \newline The condition $f(a+b) \equiv f(a)+f(b)+ab \mod h$ is equivalent to $$f(\sum_{i=1}^{n} (a_i+b_i)e_i)=f(\sum_{i=1}^{n} a_ie_i)+f(\sum_{i=1}^{n} b_ie_i)+(\sum_{i=1}^{n} a_ie_i)(\sum_{i=1}^{n} b_ie_i)$$ Using the known value of $f$, this gives $$\sum_{i=1}^{n} (a_i+b_i)f(e_i)+\binom{a_i+b_i}{2}e_i^2+\sum_{1 \leq i<j \leq n} (a_i+b_i)(a_j+b_j)e_ie_j \equiv \sum_{i=1}^{n} (a_i+b_i)f(e_i)+$$$$(\binom{a_i}{2}+\binom{b_i}{2})e_i^2+\sum_{1 \leq i<j \leq n} (a_ia_j+b_ib_j)e_ie_j+(\sum_{i=1}^{n} a_ie_i)(\sum_{i=1}^{n} b_ie_i) \mod h$$ which is an equality. Thus any choice of $f$ on a basis determines a working $f$ on all of $k_0$. As there are $2^{\text{codim}_{\mathbb{F}_2}(h)}$ choices of coset for each basis element, this gives a total of $2^{r\text{codim}_{\mathbb{F}_2}(h)}$ total choices of $f$. \end{proof}
\begin{lemma}
Let $t$ be the trace $0$ subspace of $l=\mathbb{F}_{p^k}$, i.e. the kernel of the trace map $\mathbb{F}_{p^k} \to \mathbb{F}_p$, let $l'=\mathbb{F}_{p^{k'}}$ be a subfield, and $G=\text{Gal}(l/l')$ with $g$ a generator. Set $r=|G|$ and $s=\frac{k}{r}$. Then there is an $\mathbb{F}_p$-basis for $t$ of the form $\{g^i\alpha_j, 0 \leq i \leq r-1, 1 \leq i \leq s-1\} \cup \{g^i(g-1)\alpha_s, 0 \leq i \leq r-2\}$. \end{lemma} 
\begin{proof} We start with the $\mathbb{F}_p[x]$-module structure of $l$, where $x$ acts as multiplication by $g$, which is $\prod_{i=1}^{s} \mathbb{F}_p[x]/(x^r-1)$. Note that $t$ is a $G$-invariant subspace since any $l'$-conjugate is certainly an $\mathbb{F}_p$-conjugate. Hence it also is naturally endowed with the structure of a $\mathbb{F}_p[x]$-module, and so it too has a decomposition into elementary divisors $\prod_{i=1}^{m} \mathbb{F}_p[x]/(p_i(x)^{e_i})$ according to the structure theorem (so that the $p_i$s are irreducible). \newline \newline Then note that multiplication by $x-1$ on $\prod_{i=1}^{s} \mathbb{F}_p[x]/(x^r-1)$ gives a submodule of $t$ under this isomorphism, meaning that $t$ contains the submodule $\prod_{i=1}^{s} (x-1)/(x^r-1)$ and thus a submodule isomorphic to $\prod_{i=1}^{s} \mathbb{F}_p[x]/(\frac{x^r-1}{x-1})$. Breaking both these submodules into their invariant factor decompositions and choosing any monic irreducible $p(x) \neq x-1$ in these decompositions, it follows that the dimension of the $p(x)^e$ is the same for both of these for any $e$, and so it follows that $t$ has identical $p(x)^e$-torsion, and so in particular these elementary divisors match. \newline \newline The only other possibility for a $p_i(x)$ in the decomposition for $t$ is $x-1$ itself, and by considering $(x-1)^e$ torsion in $\prod_{i=1}^{s} \mathbb{F}_p[x]/(\frac{x^r-1}{x-1})$, we see that all exponents must be at least one less than the common exponent in the decomposition of $k^+$. For dimension reasons, we must then have that all exponents are equal except for one which is one less. Thus we get an invariant factor decomposition of $t$ of the form $\prod_{i=1}^{s-1} \mathbb{F}_p[x]/(x^r-1) \oplus \mathbb{F}_p[x]/(\frac{x^r-1}{x-1})$. The element corresponding to $1$ in the last summand and $0$ elsewhere is in the kernel of the trace map $\text{Tr}_{l/l'}$, and thus is of the form $(g-1)\alpha_s$ for some $\alpha_s$ (by a counting argument or Hilbert 90), completing the proof.
\end{proof} \begin{lemma}
Let $k$ be a finite field of characteristic $2$ with $[k:\mathbb{F}_2]=n$ and choose $\alpha \in k$. Let $q=2^r$ be a prime power with $n=rm$ so that $m \geq 2$ and set $\alpha_i=\alpha^{q^i}+\alpha^{q^{i+1}}$. Let $t$ be the trace $0$ subspace. Then $\sum_{0 \leq i<j<m} \alpha_i\alpha_j \equiv \alpha^{q+1}+\alpha \mod t$ if $m$ is odd and $0 \mod t$ if $m$ is even. \begin{proof} Note that $\sum_{0 \leq i<j<m} \alpha_i\alpha_j=\sum_{0 \leq i<j<m} \alpha^{q^i+q^j}+\alpha^{q^{i+1}+q^{j+1}}+\alpha^{q^{i+1}+q^j}+\alpha^{q^i+q^{j+1}}$. First note that the pairs $\mod m$ obtained by $(i+1,j+1)$ for $0 \leq i<j<m$ are just the pairs $(i,j)$ for $1 \leq i<j \leq m$. Hence these completely overlap the $(i,j)$ pairs with $0 \leq i<j<m$ except for those with $i=0$ in the latter case and $j=m$ in the former. Hence this part of the sum just becomes $\sum_{0<i<m} \alpha^{1+q^i}+\sum_{1 \leq i<m} \alpha^{q^i+q^m}=\sum_{0<i<m} \alpha^{1+q^i}+\alpha^{q^i+q^m}=0$. \newline \newline  Hence we just need to determine 
$\sum_{0 \leq i<j<m} \alpha^{q^{i+1}+q^j}+\alpha^{q^i+q^{j+1}}$. The pairs $(i+1, j)$ obtained for $0 \leq i<j<m$ are precisely those of the form $(i,j)$ for $1 \leq i \leq j<m$ while the pairs $(i, j+1)$ obtained for $0 \leq i<j<m$ are precisely those of the form $(i,j)$ for $0 \leq i<j \leq m$ with $i+1<j$. The pairs $(i,j)$ for $1 \leq i \leq j<m$ that are not of the form $(i,j)$ for $0 \leq i<j \leq m$ with $i+1<j$ are precisely those with $i=j$ or $i+1=j$ while the pairs $(i,j)$ for $0 \leq i<j \leq m$ with $i+1<j$ that are not of the form $(i,j)$ for $1 \leq i \leq j<m$ are precisely those with $i=0$ or $j=m$. Hence all terms in the sum cancel out except these (since they all other pairs will have exactly $2$ copies), leaving $\sum_{1 \leq i<m} \alpha^{2q^i}+\sum_{1 \leq i<m-1} a^{q^i+q^{i+1}}$+$\sum_{j=2}^{m} a^{1+q^j}+\sum_{i=0}^{m-2} \alpha^{q^i+q^m}-\alpha^{1+q^m}$ (since both of the latter two sums count the case $(0,m)$). \newline \newline Overlapping the last two sums gives $\alpha^{1+q^{m-1}}+\alpha^2+\alpha^2+\alpha^{q+q^m}+\sum_{j=2}^{m-2} 2a^{1+q^j} \equiv \alpha^{1+q^{m-1}}+\alpha^{q+q^m}-\alpha^{1+q^m} \equiv \alpha^2 \mod t$. Hence the overall sum becomes $\alpha^2+\sum_{1 \leq i<m} \alpha^{2q^i}+\sum_{1 \leq i<m-1} a^{q^i+q^{i+1}}$. Note that $\alpha^{2q^i} \equiv \alpha \mod t$ for each $i$, while similarly, $\alpha^{q^i+q^{i+1}} \equiv \alpha^{q+1} \mod t$ for each $i$ (this is because the coset of $t$ is determined by the trace and taking the trace of an $\mathbb{F}_2$-conjugate gives the same result). It follows that $\alpha^2+\sum_{1 \leq i<m} \alpha^{2q^i}+\sum_{1 \leq i<m-1} \alpha^{q^i+q^{i+1}} \equiv \alpha+(m-1)\alpha+(m-2)\alpha^{q+1} \equiv m(\alpha+\alpha^{q+1}) \mod t$, as desired. \end{proof} \end{lemma}
\section{Fiber Structures}
Recall from Lemma $3.1$ that within a single Lubin-Tate extension $K_{\pi,m}$ with $m \geq n+1$, the $(-1,n)$ extensions corresponded to Galois invariant subgroups $H'$ of $\mathcal{O}_K^{\times}/(1+\mathfrak{m}^{n+1})$ such that the map $H' \to \mathcal{O}^{\times}_K/(1+\mathfrak{m}^n)$ is surjective but $1+\mathfrak{m}^{n+1} \subset H$. Thus we are looking for subgroups $H'$ of $\mathcal{O}^{\times}_K/(1+\mathfrak{m}^{n+1})$ such that the map $H' \to \mathcal{O}^{\times}_K/(1+\mathfrak{m}^n)$ is surjective. As before, by Sylow theory, we can decompose $H'$ as $\mu_{q-1} \times H$, where $H=H' \cap (1+\mathfrak{m})/(1+\mathfrak{m}^{n+1})$. The key will then be to find the $H$ such that the map $H \to (1+\mathfrak{m})/(1+\mathfrak{m}^n)$ is surjective. We will solve this problem separately in two different ranges, but first we'll prove a general theoretic lemma:
\begin{lemma}
Suppose that $G$ is a finite group and that $S \subset G$ is a subset and $H \subset G$ is a normal subgroup. Furthermore suppose that the coset map $S \to G/H$ is surjective. Given $gH \in G/H$, let the fiber at $gH$ be $S_{gH}$. Suppose that $S_{g_1H}S_{g_2H}=S_{g_1g_2H}$ and that $S_H$ is a subgroup of $G$. Then $S \subset G$ is a subgroup. Furthermore, if $S$ is a subgroup, then $S_{g_1g_2H}=S_{g_1H}S_{g_2H}$ and $S_{gH}=s_{gH}S_H$ for any fixed $s_{gH} \in S_{gH}$.
\begin{proof}
First we show that $S$ is in fact a subgroup. We just need to check closure and inverses since $S_H$ is a subgroup containing the identity. Closure follows immediately since $S_{g_1H}S_{g_2H}=S_{g_1g_2H}$. For inverses, choose $s \in S$ and suppose that $s \in S_{gH}$ for some $gH$. Then we have that $S_{gH}S_{g^{-1}H}=S_H$, and so since $e \in S_0$, we conclude that $s^{-1} \in S_{g^{-1}H}$, showing that $S$ is indeed a subgroup. \newline \newline Now suppose that $S$ is actually a subgroup. We immediately have that $S_{g_1H}S_{g_2H} \subset S_{g_1g_2H}$ since the projection map is a homomorphism. We first show that $S_{gH}=s_{gH}S_H$ for any fixed $s_{gH}$. We have that $S_{gH}S_H \subset S_{gH}$, and so $s_{gH}S_H \subset S_{gH}$, so it suffices that show that $|S_H|=|S_{gH}|$. Already we have that $|S_H| \leq |S_{gH}|$. Conversely, note that for a fixed $s_{g^{-1}H}$, we have that $S_{gH}s_{g^{-1}H} \subset S_H$, so $|S_{gH}| \leq |S_H|$, proving the claim. Then note that in general $s_{g_1H}S_{g_2H} \subset S_{g_1g_2H}$ and these sets have the same size, so $S_{g_1g_2H}=s_{g_1H}S_{g_2H} \subset S_{g_1H}S_{g_2H}$, implying that $S_{g_1g_2H}=S_{g_1H}S_{g_2H}$, as desired.
\end{proof}
\end{lemma}
Equipped with this lemma, we can now ensure that a subset of $(1+\mathfrak{m})/(1+\mathfrak{m}^{n+1})$ satisfying this compatibility property on fibers is actually equivalent to the data of a subgroup, reducing the count to counting a choice of compatible fibers. \newline \newline  We will now introduce a framework for understanding these fibers by choosing appropriate lifts. Suppose $K/\mathbb{Q}_p$ is finite with ramification index $e$, and suppose that $e \geq n$. Let $K_0$ be a finite unramified extension of $K$, and fix a uniformizer $\pi$ of $K$. Let $f=[K_0:K]$ and $k_0$ be the residue field of $K_0$. 
\begin{definition} Given some $a \in k_0$ by Hensel's Lemma there is a unique root of unity $\omega(a) \in \mathcal{O}_{K_0}$ such that its projection in the residue field is $a$, known as the Teichmüller lift. We will call these our lifts and denote them by $a'$. \end{definition}
Our lifts satisfy a crucial property in the case where $e \geq n$. \begin{lemma} Given our setup, the map $\phi: (1+Tk_0[T])/(1+T^{n+1}k_0[T]) \to (1+\mathfrak{m})/(1+\mathfrak{m}^{n+1})$ given by $1+\sum_{i=1}^{n} a_iT^i \to 1+\sum_{i=1}^{n} a'_i\pi^i$ is a Galois equivariant isomorphism. \end{lemma} \begin{proof} Galois-equivariance is immediate since the Galois action commutes with taking the Teichmüller lifts. Furthermore, every element in $(1+\mathfrak{m})/(1+\mathfrak{m}^{n+1})$ has a unique representative in the form $1+\sum_{i=1}^{n} a'_i\pi^i$, and so the map is bijective. To show it is a homomorphism, note that $$\phi(1+\sum_{i=1}^{n} a_iT^i)\phi(1+\sum_{i=1}^{n} b_iT^i)=1+\sum_{i=1}^{n} (a'_i+b'_i+\sum_{j=1}^{i-1} a'_jb'_{i-j})\pi^i$$ Then note that for each $i$, $a'_i+b'_i+\sum_{j=1}^{i-1} a'_jb'_{i-j}$ and $(a_i+b_i+\sum_{j=1}^{i-1} a_jb_{i-j})'$ are each elements in the maximal unramified subextension $M$ of $K_0/\mathbb{Q}_p$, over which $K_0$ will have ramification index $e$. Furthermore, $$a'_i+b'_i+\sum_{j=1}^{i-1} a'_jb'_{i-j}-(a_i+b_i+\sum_{j=1}^{i-1} a_jb_{i-j})'$$ lies in $\mathfrak{m}$, and thus it has valuation at least $e$. Since $e \geq n$, it follows that for each $1 \leq i \leq n$, we have $$\sum_{i=1}^{n} (a'_i+b'_i+\sum_{j=1}^{i-1} a'_jb'_{i-j})\pi^i \equiv (a_i+b_i+\sum_{j=1}^{i-1} a_jb_{i-j})'\pi^i \mod \mathfrak{m}^{n+1}$$ and so $$1+\sum_{i=1}^{n} (a'_i+b'_i+\sum_{j=1}^{i-1} a'_jb'_{i-j})\pi^i=1+\sum_{i=1}^{n} (a_i+b_i+\sum_{j=1}^{i-1} a_jb_{i-j})'\pi^i=\phi((1+\sum_{i=1}^{n} a_iT^i)(1+\sum_{i=1}^{n} b_iT^i))$$ as desired. \end{proof} We will now prove a proposition that characterizes the working fibers for our subgroup $H$ more precisely under the assumption $e \geq n$, the crucial assumption of Lemma $6.3$: \begin{prop} The data of a Galois-invariant subgroup $H$ satisfying Proposition $3.1$ for $K_0$ is equivalent to a choice of Galois-invariant subgroup $h \subset k_0^+$ and the choice of a Galois-equivariant $f: k_0^{n-1} \to k_0/h$ satisfying the relation $f(a_1,\cdots,a_{n-1})+f(a_1,\cdots,a_{n-1})+\sum_{i=1}^{n-1} a_ib_{n-i} \equiv f(c_1,\cdots,c_{n-1}) \mod h$ of Theorem $5.1$. \end{prop} \begin{proof} Firstly, note that the isomorphism of Lemma $6.3$ is both Galois-equivariant and commutes with the projection maps onto the $\mod \mathfrak{m}^n, T^n$ quotients, and so the subgroups $H$ of $(1+\mathfrak{m})/(1+\mathfrak{m}^{n+1})$ for which the natural map $H \to (1+\mathfrak{m})/(1+\mathfrak{m}^n)$ is surjective are in bijection with those of $(1+Tk_0[T])/(1+T^{n+1}k_0[T])$ satisfying the same property. \newline \newline We need to show two directions. Begin with such a function $f$ and a Galois-invariant subgroup $h \subset k_0^+$. Now we construct the fiber under $1+\sum_{i=1}^{n-1} a_iT^i$ under the projection map to be $$H_{a_1,\cdots,a_{n-1}}=1+\sum_{i=1}^{n-1} a_iT^i+f(a_1,\cdots,a_{n-1})T^n$$ where $f(a_1,\cdots,a_{n-1})$ corresponds to some coset of $h$ in $k_0$. To verify the hypotheses of Lemma $6.1$, it suffices to show that $H_0$ is a subgroup and that $H_{a_1,\cdots,a_{n-1}}H_{b_1,\cdots,b_{n-1}}=H_{c_1,\cdots,c_{n_1}}$. The former follows since $f(0)=0$, which implies that $H_0=1+T^nh$, which is in fact of subgroup of $(1+Tk_0[T])/(1+T^{n+1}k_0[T])$, while $$H_{a_1,\cdots,a_{n-1}}H_{b_1,\cdots,b_{n-1}}=(1+\sum_{i=1}^{n-1} a_iT^i+f(a_1,\cdots,a_{n-1})T^n)(1+\sum_{i=1}^{n-1} b_iT^i+f(b_1,\cdots,b_{n-1})T^n)=$$$$1+\sum_{i=1}^{n-1} c_iT^i+(\sum_{i=1}^{n-1} a_ib_{n-i}+f(a_1,\cdots,a_{n-1})+f(b_1,\cdots,b_{n-1}))T^n=1+\sum_{i=1}^{n-1} c_iT^i+f(c_1,\cdots,c_n)T^n$$ as desired, showing that this $H$ indeed yields a subgroup. To show that the subgroup is Galois-invariant, simply note that since $f$ is Galois-equivariant, $\sigma(H_{a_1,\cdots,a_{n-1}})=1+\sum_{i=1}^{n-1} \sigma(a_i)T^i+\sigma(f(a_1,\cdots,a_{n-1}))T^n=1+\sum_{i=1}^{n-1} \sigma(a_i)T^i+f(\sigma(a_1),\cdots,\sigma(a_{n-1}))T^n=H_{\sigma(a_1),\cdots,\sigma(a_{n-1})}$ showing Galois invariance. \newline \newline Conversely, suppose we have such a subgroup $H$, viewed via its corresponding subgroup $H' \subset (1+Tk_0[T])/(1+T^{n+1}k_0[T])$ and recall that by Lemma $6.1$ that $H'$ is characterized by its fibers $H'_{a_1,\cdots,a_{n-1}}$ as long as $H'_0$ is a subgroup and $H'_{a_1,\cdots,a_{n-1}}H'_{b_1,\cdots,b_{n-1}}=H'_{c_1,\cdots,c_{n-1}}$. Note that $H'_0$ will correspond to a subgroup of $(1+T^nk_0[T])/(1+T^{n+1}k_0[T]) \cong k_0^+$ by definition, and so we can identify it with the subgroup $h$ of $k_0^+$ and thus write it as $H_0=1+T^nh$. Note that $H'_0$ consists of the cosets with representatives $x$ such that $v_T(x-1) \geq n$, and since $\sigma$ preserves valuations and $H'$ is Galois invariant, it follows that $\sigma(H'_0)=H'_0$ and thus $h$ is Galois-invariant. By Lemma $6.1$, any other fiber $H'_{a_1,\cdots,a_{n-1}}$ is of the form $h'_{a_1,\cdots,a_{n-1}}H'_0$ for any fixed $h'_{a_1,\cdots,a_{n-1}} \in H'_{a_1,\cdots,a_{n-1}}$. We can express any given element $h'_{a_1,\cdots,a_{n-1}}$ in the form $1+\sum_{i=1}^{n-1} a_iT^i+xT^n$ for some $x$, and then we can take the coset of $x+h$ in $k_0$, which we call $f(a_1,\cdots,a_{n-1})$. It immediately follows that $$H'_{a_1,\cdots,a_{n-1}}=h_{a_1,\cdots,a_{n-1}}H'_0=1+\sum_{i=1}^{n-1} a_iT^i+f(a_1,\cdots,a_{n-1})T^n$$ 
and so we can think of the fibers as being of the form $$H'_{a_1,\cdots,a_{n-1}}=1+\sum_{i=1}^{n-1} a_iT^i+f(a_1,\cdots,a_{n-1})T^n$$ Then since $H$ is a subgroup, we will have that $$1+\sum_{i=1}^{n-1} c_iT^i+f(c_1,\cdots,c_{n-1})T^n=H'_{c_1,\cdots,c_{n-1}}=H'_{a_1,\cdots,a_{n-1}}H'_{b_1,\cdots,b_{n-1}}=$$$$1+\sum_{i=1}^{n-1} c_iT^i+(\sum_{i=1}^{n-1} a_ib_{n-i}+f(a_1,\cdots,a_{n-1})+f(b_1,\cdots,b_{n-1}))T^n$$ and thus we conclude that $$f(c_1,\cdots,c_{n-1}) \equiv \sum_{i=1}^{n-1} a_ib_{n-i}+f(a_1,\cdots,a_{n-1})+f(b_1,\cdots,b_{n-1}) \mod h$$ It follows that $H'_{c_1,\cdots,c_{n-1}}=H'_{a_1,\cdots,a_{n-1}}H'_{b_1,\cdots,b_{n-1}}$ is equivalent to $f(c_1,\cdots,c_{n-1}) \equiv \sum_{i=1}^{n-1} a_ib_{n-i}+f(a_1,\cdots,a_{n-1})+f(b_1,\cdots,b_{n-1}) \mod h$ which defines our group structure. \end{proof}
The reason why this lemma is so important is that it ultimately, utilizing the choice of lifts in the ramified setting, collapses a nonlinear problem that lives in $(1+\mathfrak{m})/(1+\mathfrak{m}^{n+1})$ into a problem that inherently lives in the vector space $k_0^+$, allowing us to utilize the tools of Section $5$ in our counting problem.
\section{Main Results} 
\begin{theorem}
Suppose $K/\mathbb{Q}_p$ is finite with ramification index $e$, inertia degree $f'$, and residue field $k$. Let $f$ be a positive integer and let $q=p^{f'}, q'=q^f$, $k_0=\mathbb{F}_{q'}$, and $G=\text{Gal}(k_0/k)$. Then let $S$ be the set of $G$-invariant subspaces of $k_0$ and $t \in S$ to be the trace $0$ subspace under $Tr_{k_0/k}$. Let $b=n-1-\lfloor \frac{n-1}{p} \rfloor$. Then the number of $(-1,n)$ Galois extensions $L/K$ of inertia degree $f$ is $$\begin{cases} 0 & p|n, e \geq n \\ q^{bf'f+1}\sum_{h \in S} \frac{|h \cap t|}{|h|^{bf'+1}} & e \geq n, p \nmid n \\ 0 & e<n, p>e+1 \end{cases}$$ \end{theorem}
\begin{proof} Let $K_0/K$ be the maximal unramified subextension of $L/K$. First note that $G=\text{Gal}(L/K_0)$ has single wild ramification jump at $n$, which means that the extension is totally wildly ramified. Thus by Proposition $2.1$, we see that $G_n/G_{n+1} \cong \text{Gal}(L/K_0)$ is the direct sum of cyclic groups of order $p$. In particular, $G$ is abelian, and so by class field theory we may attach a norm group $\text{Nm}_{L/K_0}(L^\times)$ to it. As $L/K$ is totally ramified, we may let $\pi \in \text{Nm}_{L/K_0}(L^\times)$ be a uniformizer of $\mathcal{O}_{K_0}$. It then follows that $\text{Nm}_{L/K_0}(L^\times)=\pi^{\mathbb{Z}}\text{Nm}_{L/K_0}(\mathcal{O}_L^\times)$. Let $H=\text{Nm}_{L/K_0}(\mathcal{O}_L^\times)$. By Proposition $2.2$, since $G$ is abelian, the extension $L/K$ being Galois is equivalent to having the norm group corresponding to the abelian extension $L/K$ under class field theory to be fixed by $\text{Gal}(K_0/K)$. This means that $\pi^{\mathbb{Z}}H$ is invariant under the Galois action. \newline \newline Since $\sigma$ preserves valuations, $\pi^{\mathbb{Z}}H$ is invariant under $\text{Gal}(K_0/K)$ iff $\pi^fH$ is for each integer $f$. As $K_0/K$ is unramified, let $\pi=\pi'u$ for some uniformizer $\pi'$ of $\mathcal{O}_{K}$ and $u \in \mathcal{O}_{K_0}^{\times}$. The group $\mathcal{O}_{K_0}^\times/H$ has finite order, and so choosing $f=|\mathcal{O}_{K_0}^\times/H|$, Lagrange implies that $\pi^fH=\pi'^fH$. Thus to be Galois invariant in this case just means that $H$ is Galois invariant. Now knowing that $H$ is Galois invariant, we see that $\pi^{\mathbb{Z}}H$ is invariant precisely if $\sigma(u)/u \in H$ for each $\sigma \in G=\text{Gal}(K_0/K)$. Thus for a given Galois invariant $H$, it suffices to find the number of classes $u \in \mathcal{O}_{K_0}^\times/H$ that are also Galois invariant. As $H$ is Galois invariant, $\mathcal{O}_{K_0}^\times/H$ naturally obtains the structure of a $G$-module. 
\newline \newline By Proposition $6.4$, the subgroups we seek to count that are contained within a particular Lubin-Tate extension are characterized by a choice of Galois-invariant subgroup $h \subset k_0^+$ and a Galois-equivariant function $f: k_0^{n-1} \to k_0/h$ such that $f(a_1,\cdots,a_{n-1})+f(a_1,\cdots,a_{n-1})+\sum_{i=1}^{n-1} a_ib_{n-i} \equiv f(c_1,\cdots,c_{n-1}) \mod h$. By Theorem $5.1$, the number of functions satisfying these conditions $q^{b\text{codim}_{\mathbb{F}_p}(h)}$ if $p \nmid n$ and $0$ if $p|n$. \newline \newline Now since the projection map $H \rightarrow \mathcal{O}_{K_0}^\times/(1+\mathfrak{m}^n)$ is surjective, any element of $\mathcal{O}_{K_0}^\times/H$ has a coset representative of the form $1+\pi'^nx$. We want to compute the number of $G$-invariant points of $\mathcal{O}_{K_0}^\times/H$ given a choice of $H$. Recalling that $\pi'$ is a uniformizer of $\mathcal{O}_K$ preserved by the Galois action, we may view $x$ as an element of $k$ since shifting $x$ by something in $\mathfrak{m}$ does not change its coset. Such an element $x$ is then precisely defined by its coset in $k_0/h$. The action of $G$ on $k_0/h$ restricts to the action of $G_1=\text{Gal}(k_0/k)$ on $k_0/h$, and so we just seek the number of $G_1$ invariant fixed points of $k_0/h$ for a given choice of $h$. \newline \newline Now note that its coset $x+h$ is invariant under $G_1$ iff it is invariant under a generator $\sigma$, meaning that we just need $x^q-x \in h$. Thus we seek the number of elements $x+h$ of $k_0/h$ such that $x^q-x \equiv 0 \mod h$, where $q=|k|$. The map $x \rightarrow x^q-x$ is a linear map $k_0 \rightarrow k_0$ with kernel consisting of the elements of $\mathbb{F}_q$ and image $t$, where $t$ is the trace $0$ subspace (i.e. the kernel of $\text{Tr}_{k_0/k}$), since anything in the image is in the kernel of the trace map and $t$ and the image have the same order. The subspace of $h$ in the image of this map is then $h \cap t$. Each of these images is attained $q$ times, so the total number of images in $h$ is $q|h \cap t|$. The total number of cosets of $h$ is then $q|h \cap t|/|h|$, and so this is the number of fixed points. For each $|h|$, there are $q^{b\text{codim}_{\mathbb{F}_p}(h)}=(\frac{|k_0|}{|h|})^{bf'}$ ways to extend it to a $G$-invariant subgroup $H$ of $\mathcal{O}_{K_0}^\times$ that is surjective on the projection map, and so to get the total number of extensions, we sum this over all possible $|h|$. Thus the answer is general is $$q\sum_{h \subset k_0^+} (\frac{|k_0|}{|h|})^{bf'}\frac{|h \cap t|}{|h|}=q^{bff'+1}\sum_{h \subset k_0^+} \frac{|h \cap t|}{|h|^{bf'+1}}$$ where $h$ ranges over all $\text{Gal}(k_0/k)$ stable subspaces of $k_0^+$, as desired. \newline \newline Now it remains to handle the case when $e<n, p>e+1$. By Proposition $2.1$, $G_n \cong G_n/G_{n+1} \cong \text{Gal}(L/K_0)$ is the direct sum of cyclic groups of order $p$. By Artin reciprocity and noting that $L/K_0$ is totally ramified, we have that $\text{Gal}(L/K_0) \cong \mathcal{O}^\times_{K_0}/\text{Nm}_{L/K_0}(\mathcal{O}_L^\times)$. This means that $(\mathcal{O}_{K_0}^\times)^p \subset \text{Nm}_{L/K_0}(\mathcal{O}_L^\times)$ However, by Remark $2.5$, we have that $(1+\mathfrak{m})^p \supset 1+\mathfrak{m}^{e+1} \supset 1+\mathfrak{m}^n$. Combined with the surjectivity of $H \to (1+\mathfrak{m})/(1+\mathfrak{m}^n)$, this implies that $H$ is not a proper subgroup of $1+\mathfrak{m}$, contradicting Lemma $3.1$, as desired. \end{proof}
In Theorem $7.1$, we crucially assumed that $e \geq n$ in order for our coset machinery to work properly. In fact, the count depends entirely on the content of the residue field extension $k_0/k$ of the maximal unramified subextension, and so it depends entirely on the residue field extension. 
\begin{corollary} Suppose that $K/\mathbb{Q}_p$ is totally ramified so that $f'=1$. For a positive integer $f$, set $m(x)=\frac{x^f-1}{(x-1)^{p^{v_p(f)}}} \in \mathbb{F}_p[x]$, $\text{deg}(m)=f-p^{v_p(f)}=d$, and $\zeta_m(s)=\sum_{i=0}^{d} \frac{a_i}{p^{is}}$ be the $\zeta$-function for the ring $S=\frac{\mathbb{F}_p[x]}{(m)}$ (equivalently, $a_i$ is the number of degree $i$ monic factors of $m$). Then the number of $(-1,n)$ Galois extensions $L/K$ with inertia degree $f$ is $$\begin{cases} 0 & p|n, e \geq n \\ \frac{(p^{(f+1)b}-p^{fb}+p^{(f-1)b+1}-p^{db+1})}{p^b-1}\zeta_m(b)  & e \geq n, p \nmid n \\ 0 & e<n, p>e+1 \end{cases}$$ \end{corollary} 
\begin{proof} Only the case $e \geq n, p \nmid n$ needs to be shown. By Theorem $7.1$ and upon specializing to $f'=1$ and $q=p^{f'}=p$, the number of such extensions is $$p^{bf+1}\sum_{h \subset k_0^+} \frac{|h \cap t|}{|h|^{bf+1}}$$ Now fix $|h|$. We want to determine $\sum |h \cap t|$ over all subspaces $h$ of $k_0$ with $|h|$ of a given size. For this, we use the Galois module structure. \newline \newline First suppose that $K/\mathbb{Q}_p$ is totally ramified. Since $K/\mathbb{Q}_p$ is totally ramified, $\mathbb{F}_{p^f}$ naturally has a structure as a $\mathbb{F}_p[x]$-module, decomposing as $\mathbb{F}_p[x]/(x^f-1)$. The virtue of this is that $\mathbb{F}_p[x]$-submodules, i.e. Galois invariant $\mathbb{F}_p$-subspaces, correspond precisely to ideals of $R=\mathbb{F}_p[x]/(x^f-1)$. \newline \newline Under this correspondence, the subspace $t'$ is just the ideal $(x-1)R$, and so $|h \cap t'|=|(x-1) \cap I|$. But now $I$s is necessarily a principal ideal corresponding to a monic factor $f$ of $x^n-1$ in $\mathbb{F}_p[x]$. Let $I=(f(x))$. If $x-1|f(x)$, then $I \subset (x-1)$ and so $|(x-1) \cap I|=|I|$. Otherwise, $(x-1) \cap I=(f(x)(x-1))$ and so $|(x-1) \cap I|=|I|/p$. But then $|h|=|I|$ and so the sum we seek is $|h|\sum_{i=0}^{1} \frac{N_{i,t}}{p^i}$. Let $|h|=p^j$. Thus the total sum over $h$ with $|h|=p^j$ becomes $p(\sum_{i=0}^{1} \frac{N_{i,j}}{p^i})$. \newline \newline Now this is the total number of fixed points for a given choice of $|h|$. Breaking the original sum into parts by $|h|$, we get $$p^{fb+1} \sum_{j=0}^{n-1} \sum_{i=0}^{1} \frac{N_{i,j}}{p^{i+jb}}$$ Now we determine $N_{i,j}$ more explicitly. Firstly, if $i=1$, then $N_{1,j}$ counts the number of degree $t$ factors indivisible by $x-1$, which is $a_{d-t}=a_j$. If $i=0$, then we now restrict to those divisible by $x-1$, which is $a_{j-1}+a_{j-2}+\cdots+a_{j-p^{v_p(f)}}$. Thus the part with $i=1$ becomes $p^{fb}\sum_{j=0}^{d} \frac{a_j}{p^{jb}}$ while the part with $i=0$ is just $p^{fb+1}(\sum_{j=0}^{d}\frac{a_t}{p^{jb}})(\sum_{j=1}^{p^{v_p(f)}} \frac{1}{p^{jb}})$. Combining the two, we get $\frac{(p^{f+1)b}-p^{fb}+p^{(f-1)b+1}-p^{db+1})}{p^b-1}\zeta_m(b)$, as desired. \end{proof}
Note that the only cases not handled by Theorem $7.1$ are those where $e<n, p-1 \leq e$, meaning the chain of inequalities $p-1 \leq e<n$ must hold. This leaves only finitely many $(e,p)$ pairs for a given $n$, effectively closing the entire range of the problem. When $n=2$, the only exceptional case is $e=1,p=2$, which we will now settle: 
\begin{theorem}
Suppose $K/\mathbb{Q}_2$ is unramified of degree $n$. Then the number of $(-1,2)$ Galois extensions $L/K$ with inertia degree $r$ is $2^{n+1}$ (independent of $r$). \end{theorem}
\vspace{-0.6\baselineskip}
\begin{proof} Let $K_0/K$ be the maximal unramified subextension with $K_0/K$ having residue field extension $k_0/k$. As before, we have that $(\mathcal{O}_{K_0}^\times)^2 \subset \text{Nm}_{L/K_0}(\mathcal{O}_L^\times)$. We will determine precisely $(\mathcal{O}_K^\times)^2$. Note that $(1+2z)^2=1+4(z^2+z)$. If we can find $x$ such that $z^2+z=y$, then we have that $1+2z$ is a solution to $x^2-(1+4y)=0$, and so Hensel's Lemma implies that we can find $w$ with $w^2=1+4y$. Thus the squares in $1+2\mathcal{O}_K$ are precisely those of the form $1+4(z^2+z) \mod 8$. In particular, we deduce that if $h$ is the image of $H \cap (1+4\mathcal{O}_{K_0}) \rightarrow (1+4\mathcal{O}_{K_0})/(1+8\mathcal{O}_{K_0}) \cong k_0^+$, then $h$ contains the subspace consisting of all values $z^2+z, z \in k$. The map $z \to z^2+z$ is a homomorphism on $k$ with kernel $0,1$, and thus the image has size $\frac{|k_0|}{2}$, meaning that this subspace is an index $2$ subgroup of $k_0^+$. In particular, we must have that $h$ is precisely this subspace or else $h$ would be the whole of $k_0^+$ which would imply as before that $H=\mathcal{O}_{K_0}^\times$, which is again not proper, contradicting Lemma $3.1$. Thus $h=t$, the trace $0$ subspace.
\vspace{0.6\baselineskip}

\noindent 
By Lemma $4.3$, we can choose an $\mathbb{F}_2$-basis for $t$ of the form $\{g^i\alpha_j, 0 \leq i \leq r-1, 1 \leq j \leq n-1\} \cup \{g^i(g-1)\alpha_n, 0 \leq i \leq r-2\}$, where $g$ is a generator of $\text{Gal}(k_0/k)$. Denote these elements as $e_1,\cdots,e_{rn-1}$. Lift these basis elements to roots of unity $\omega_1,\cdots,\omega_{rn-1}$ in $\mathcal{O}_K^{\times}$. By Lemma $4.3$, we can choose an $\mathbb{F}_2$-basis for $t$ of the form $\{g^i\alpha_j, 0 \leq i \leq r-1, 1 \leq j \leq n-1\} \cup \{g^i(g-1)\alpha_n, 0 \leq i \leq r-2\}$, where $g$ is a generator of $\text{Gal}(k/k_1)$. Denote these elements as $e_1,\cdots,e_{rn-1}$. Lift these basis elements to roots of unity $\omega_1,\cdots,\omega_{rn-1}$ in $\mathcal{O}_K^{\times}$. Each $a \in k$ can uniquely written as $\sum z_ie_i$, where $z_i \in \mathbb{F}_2$, and so there is a unique lift of $a$ of the form $\sum z_i\omega_i$, where $0 \leq z_i \leq 1$. Let $a'=\sum a_i\omega_i$ be this lift of $a$. \newline \newline Now set $H_a=1+2a'+4f(a)$ for some coset $f(a)$ of $t$ in $k$. Getting a subgroup structure is equivalent to having $H_aH_b=H_{a+b}$ for all $a,b \in k$. \newline \newline First suppose that $a,b \in t$. The lift $(a+b)'$ for $a+b$ will differ from $a'+b'$ by $2(\sum e_i')$, where the sum ranges over some subset of $e_1',\cdots,e_{rn-1}'$. Thus $1+2(a'+b')+4ab+4f(a)+4f(b)=1+2(a+b)'+4(\sum e_i')+4ab+4f(a)+4f(b)$. The key point is that now $\sum e_i'$ will be an element of $h$ and thus will not change the coset of $h$ dictated by the values of $f$, and thus $1+2(a+b)'+4f(a+b)=1+2(a'+b')+4ab+4f(a)+4f(b)=1+2(a+b)'+4(\sum e_i')+4ab+4f(a)+4f(b)=1+2(a+b)'+4ab+4f(a)+4f(b)$, and so we deduce that $f(a+b) \equiv f(a)+f(b)+ab \mod t$ for all $a,b \in t$.
\newline \newline 
Case 1: $r$ is odd. \newline \newline In this case, take some element $\beta$ of $k$ such that $\text{Tr}_{k/\mathbb{F}_2}(\beta)=1$ and let $e_{rn}=\beta$. Then note that $\text{Tr}_{k_0/\mathbb{F}_2}(\beta)=[K:K_1]\text{Tr}_{k/\mathbb{F}_2}(\beta)=1$, and so $\beta \not \in t$, showing that $\beta$ completes the basis. Then if both $a,b \not \in t$, the element $(a+b)'-(a'+b')$ will now be of the form $2(\sum e_i')$, where it includes $e_{rn}'$, and so now $\sum e_i'$ will indeed change the coset. Hence we deduce that $f(a+b) \equiv f(a)+f(b)+ab+e_{rn} \mod h$. Let $l$ be the subspace of $k$ spanned by $e_1,\cdots,e_{r(n-1)}$. Applying Lemma $4.1$ upon noting that $l$ has a Galois  invariant basis and stays in $t$ upon squaring (since $t$ maps to itself upon squaring), meaning that $l \sim t$, we conclude that there are $2^{n-1}$ Galois equivariant functions satisfying $f(a+b) \equiv f(a)+f(b)+ab \mod t$ on $l$. \newline \newline We will now show that each of these equivariant functions has exactly two extensions to $k_0^+$. First note that $f((g-1)\alpha_n)=f(\alpha_n^q+\alpha_n)$. In order for the functional equation to be satisfied and to have Galois equivariance, we must have that $f(\alpha_n^q+\alpha_n) \equiv 2f(\alpha_n)+\alpha_n^{q+1} \equiv \alpha_n^{q+1} \mod h$ if $\alpha_n \in t$ and $f(\alpha_n^q+\alpha_n) \equiv 2f(\alpha_n)+\alpha_n^{q+1}+\alpha_n \equiv \alpha_n^{q+1}+\alpha_n \mod h$ if $\alpha_n \not \in t$. Note that the former is just $\alpha_n^{q+1}+\alpha_n$ since $\alpha_n \in t$, so we have that $f(\alpha_n^q+\alpha_n)=\alpha_n^{q+1}+\alpha_n$ in all cases. Then since the Galois group preserves $t$, Galois equivariance forces $f(\alpha_n^{q^i}+\alpha_n^{q^{i-1}}) \equiv \alpha_n^{q+1}+\alpha_n \mod h$ for each $i$. Finally, we let $f(\beta)$ be either of the two possibilities. This defines $f$ on all of $k^+$ via the functional equation, and by the proof of Lemma $4.1$, will give a well-defined function. \newline \newline Hence we just need to check that each such function is Galois equivariant. Write $a=a_{rn}\beta+\sum_{i=1}^{n(r-1)} a_ie_i+c$ and let $b=a_{rn}\beta+\sum_{i=1}^{n(r-1)} a_ie_i$. Note that $\sigma(f(a))=\sigma(f(b))+\sigma(f(c))+\sigma(bc)$. Since the basis elements in the expansion of $b$ form a Galois invariant set, the proof of Lemma $4.2$ implies that $\sigma(f(b))=f(\sigma(b))$. Since we know that $f$ satisfies the functional equation, we must have that $f(\sigma(a))=f(\sigma(b))+f(\sigma(c))+\sigma(b)\sigma(c)$. Hence it suffices to show that $\sigma(f(c))=f(\sigma(c))$. \newline \newline By definition and Galois equivariance on the basis, we have $c=\sum_{i=r(n-1)}^{rn-1} a_ie_i=\sum_{i=0}^{r-2} a'_ig^i(g-1)\alpha_n$ and so $$g(f(c))=\sum_{i=0}^{r-2} a'_ig(f(g^i(g-1)\alpha_n))+\sum_{0 \leq i<j \leq r} g(a'_ia'_j(g^i(g-1)\alpha_n)(g^j(g-1)\alpha_n))=$$ $$\sum_{i=0}^{r-2} a'_i(f((g-1)\alpha_n))+\sum_{0 \leq i<j \leq r-2} a'_ia'_j(g^{i+1}(g-1)\alpha_n)(g^{j+1}(g-1)\alpha_n))$$ Next note that $gc=a'_{r-2}(g-1)\alpha_n+\sum_{i=0}^{r-3} (a'_i+a'_{r-2})g^{i+1}(g-1)\alpha_n$ and so $$f(gc)=f(\sum_{i=0}^{r-3} a'_ig^{i+1}(g-1)\alpha_n)+a'_{r-2}f(\sum_{i=0}^{r-2} g^i(g-1)\alpha_n)+$$$$a'_{r-2}(\sum_{i=0}^{r-3} a'_ig^{i+1}(g-1)\alpha_n)(\sum_{i=0}^{r-2} g^i(g-1)\alpha_n)$$ Now we know that $$f(\sum_{i=0}^{r-3} a'_ig^{i+1}(g-1)\alpha_n)=(\sum_{i=0}^{r-3} a'_if(g^{i+1}(g-1)\alpha_n))+$$$$\sum_{0 \leq i<j \leq r-3} a'_ia'_j(g^{i+1}(g-1)\alpha_n)(g^{j+1}(g-1)\alpha_n))$$ Combining everything, it remains to show that $f(\sum_{i=0}^{r-2} g^i(g-1)\alpha_n) \equiv f((g-1)\alpha_n) \mod h$. Expanding out the left hand side gives $(r-1)f((g-1)\alpha_n)+\sum_{0 \leq i<j \leq r-2} g^i(g-1)\alpha_n(g^j(g-1)\alpha_n)$, and so since $f((g-1)\alpha_n) \equiv \alpha^{q+1}+\alpha \mod h$, it remains to show that $\sum_{0 \leq i<j \leq r-2} g^i(g-1)\alpha_n(g^j(g-1)\alpha_n) \equiv r(\alpha^{q+1}+\alpha) \mod h$. By Lemma $4.4$, we have that $\sum_{0 \leq i<j \leq r-1} g^i(g-1)\alpha_n(g^j(g-1)\alpha_n)=r(\alpha_n^{q+1}+\alpha_n)$. But then rewriting the sum gives $$\sum_{0 \leq i<j \leq r-2} g^i(g-1)\alpha_n(g^j(g-1)\alpha_n)+(\sum_{i=0}^{r-2} g^i(g-1)\alpha_n)(g^{r-1}(g-1)\alpha_n)=$$$$\sum_{0 \leq i<j \leq r-2} g^i(g-1)\alpha_n(g^j(g-1)\alpha_n)+((1+g^{r-1})\alpha_n)(g^r-g^{r-1})\alpha_n$$ However, $((1+g^{r-1})\alpha_n)(g^r-g^{r-1})\alpha_n \equiv (\alpha+\alpha^{q^{r-1}})(\alpha^{q^r}+\alpha^{q^{r-1}}) \equiv (\alpha+\alpha^{q^{r-1}})^2=\alpha+\alpha^{q^{r-1}} \equiv 0 \mod h$, and so we conclude that $\sum_{0 \leq i<j \leq r-2} g^i(g-1)\alpha_n(g^j(g-1)\alpha_n) \equiv r(\alpha_n^{q+1}+\alpha_n) \mod h$, completing the proof. \newline \newline
Case 2: $r$ is even. \newline \newline Consider our basis $\{g^i\alpha_j, 0 \leq i \leq r-1, 1 \leq j \leq n-1\} \cup \{g^i(g-1)\alpha_n, 0 \leq i \leq r-2\}$ for $t$ coming from Lemma $4.3$. This corresponds to the $\mathbb{F}_2[x]$-module decomposition $t \cong \bigoplus_{i=1}^{n-1} \mathbb{F}_2[x]/(x^r-1) \oplus \mathbb{F}_2[x]/(x^{r-1}+\cdots+1)$. We claim that $\alpha_n \not \in t$. Suppose that $\alpha_n \in t$. Then we can write $\alpha_n=(g-1)b$ for some $b$, and so then $(g-1)\alpha_n=(g-1)^2b$. Since $r$ is even, we know that $x-1|x^{r-1}+\cdots+1$, and so it follows that $(g-1)^2b$ is annihilated by the element $\frac{g^{r-1}+\cdots+1}{g-1}$, contradiction since $(g-1)\alpha_n$ corresponds to $1$ in the last component under the $\mathbb{F}_2[x]$-module isomorphism. \newline \newline Now using the exact same reasoning as in Case $1$, we can define our function $f$ on the subspace $l$ corresponding to the Galois equivariant part of the basis for $t$ and note that that $f((g-1)\alpha_n) \equiv f(\alpha_n^q+\alpha_n) \equiv \alpha_n^{q+1}+\alpha_n \mod h$, this time using that $\alpha_n \not \in t$. Furthermore, since $\alpha_n \not \in t$, we can use it to complete our basis for $k^+$. There are then $2$ choices for $f(\alpha_n)$, each defining $f$ on all of $k_0^+$ via the functional equation, and by the proof of Lemma $4.1$ implies that we get a well-defined function in each case. 
\newline \newline Hence we just need to check that each function is Galois equivariant. The proof of equivariance in the case when the $\alpha_n$ coefficient is $0$ is identical to the proof in Case $1$. Hence we just need to show equivariance in the case where the $\alpha_n$ coefficient is $1$. Write $a=a'+\alpha$ so that $a' \in t$. Noting that we already have equivariance on $t$, we have that $\sigma(f(a)) \equiv \sigma(f(a'))+\sigma(f(\alpha))+\sigma(a'\alpha) \equiv f(\sigma(a')+\sigma(a')\sigma(\alpha)+\sigma(f(\alpha)) \mod h$. 
Hence we just need to show that $\sigma(f(\alpha)) \equiv f(\sigma(\alpha)) \mod h$. Note that $$f(g(\alpha)) \equiv f(\alpha+(g-1)\alpha)) \equiv f(\alpha)+f((g-1)\alpha)+\alpha(g-1)\alpha \equiv$$$$ f(\alpha)+\alpha^{q+1}+\alpha+\alpha(\alpha^q+\alpha) \equiv f(\alpha)+\alpha^2+\alpha \equiv f(\alpha) \mod h$$ completing the proof. \newline \newline Hence in each case, we have $2^n$ Galois equivariant functions for $h=t$, which shows that there $2^n$ possibilities for $h$ in each case. Now since the quotient $k_0^+/t$ consists of two cosets and each is fixed (since $t$ itself is), we conclude that there are $2$ fixed points. Hence following the proof of Theorem $5.1$, we conclude that there are precisely $2^{n+1}$ extensions, as desired. \end{proof}
\section*{Acknowledgments}
The author would like to thank Professor Matthias Flach for reading over a draft and leaving some comments.

\end{document}